\documentclass[10pt,a4paper]{article} 
\title{A special value of the spectral zeta function \\
of the non-commutative harmonic osciallators }
\author{Hiroyuki Ochiai
\thanks{The research of the author is supported in part by a Grant-in-Aid for 
Scientific Research (B) 15340005 from the Ministry of Education,
Culture, Sports, Science and Technology.\newline
Mathematics Subject Classification; Primary 11M36, Secondary 33C20, 33C75. \newline
Abbreviated title: A special value of spectral zeta
} 
} 
\usepackage{amsmath}
\usepackage{amssymb}
\newcommand{\R}{\mathbf{R}}

\newlength{\theoremindent}
\newtheorem{theorem}{\hspace*{\theoremindent}{\bf Theorem}}
\newtheorem{proposition}[theorem]{\hspace*{\theoremindent}{\bf Proposition}}

\newtheorem{lemma}[theorem]{\hspace*{\theoremindent}{\bf Lemma}}

\newtheorem{remark}[theorem]{\hspace*{\theoremindent}{\bf Remark}}

\date{}

\newcommand{\pz}{\partial_z}
\newcommand{\pt}{\partial_t}

\begin{document}
\maketitle


\noindent{\bf Abstract}:
The non-commutative harmonic oscillator is
a $2\times2$-system of harmonic oscillators
with a non-trivial correlation.
We write down explicitly
the special value at $s=2$ of the spectral zeta function
of the non-commutative harmonic oscillator
in terms of the complete elliptic integral of the first kind,
which is a special case of a hypergeometric function.

\section{Introduction}

The non-commutative harmonic oscillator $Q=Q(x, \partial_x)$
is defined to be the second-order ordinary differential operator
\[
Q(x,\partial_x) = 
\left[\begin{array}{cc}
\alpha & 0 \\
0 & \beta 
\end{array}\right]
(-\frac{\partial_x^2}{2} + \frac{x^2}{2} ) +
\left[\begin{array}{cc}
0 & -1 \\
1 & 0
\end{array}\right]
(x \partial_x + \frac12).
\]
The first term is 
two harmonic oscillators, which are mutually independent,
with the scaling constant $\alpha>0$ and $\beta>0$,
while the second term is considered to be
the correlation with a self-adjoint manner.
The spectral problem is a $2 \times 2$ system of the
ordinary differential equations
\[
Q(x,\partial_x) u(x) = \lambda u(x)
\]
with an eigenstate $u(x) = \left[\begin{array}{c}
u_1(x) \\
u_2(x) 
\end{array}\right] 
\in L^2(\R)^{\oplus 2}$
and a spetrum $\lambda \in \R$.
It is known \cite{PW1} that
under the natural assumption $\alpha\beta>1$
on the positivity,
which is also assumed in this paper,
the operator $Q$ defines a positive,
self-adjoint operator 
with a discrete spectum
\[
(0 <) \lambda_1 \le \lambda_2 \le \cdots \to +\infty
\]
The corresponding spectral zeta function
is defined to be 
\[
\zeta_Q(s) = \sum_{n=1}^\infty \lambda_n^{-s}.
\]
An expression of the special value $\zeta_Q(2)$
is obtained in \cite{IW}
in terms of a certain contour integral
using the solution of a singly confluent type Heun differential equation.
It would be indicated that these special values are complicated enough
and highly transcendental as reflecting the 
transcendence of the spectra of the non-commutative harmonic oscillator.

However, in this paper, we prove the following simple expression:
\begin{equation}
\label{eq:main1}
\zeta_Q(2) 
= \frac{\pi^2}{4} \frac{(\alpha^{-1}+\beta^{-1})^2}{(1-\alpha^{-1}\beta^{-1})}
\left( 1+ \left(\frac{\alpha^{-1}-\beta^{-1}}{\alpha^{-1}+\beta^{-1}}
{}_2F_1\left(\frac14,\frac34;1; \frac{1}{1-\alpha\beta} \right) \right)^2 \right).
\end{equation}
where $_2F_1$ is the  Gauss hypergeometric series.
We also derive the following representation which involves 
the complete elliptic integral of the first kind  as
\begin{equation}
\label{eq:main2}
\zeta_Q(2) 
= \frac{\pi^2}{4} \frac{(\alpha^{-1}+\beta^{-1})^2}{(1-\alpha^{-1}\beta^{-1})}
\left( 1+ \left(\frac{\alpha^{-1}-\beta^{-1}}{\alpha^{-1}+\beta^{-1}}
\int_0^{2\pi} \frac{d\theta}{2\pi\sqrt{1+ (\cos\theta)/\sqrt{1-\alpha\beta} }}\right)^2  \right).
\end{equation}
In this sense, the special value $\zeta_Q(2)$
is written in terms of a hypergeometric series,
which is more tractable and many of its properties are known.
Note that each spectrum is related 
with the monodromy problem of 
Heun's differential equation,
which is far from hypergeometric,
see \cite{O}, \cite{O2}.
Only the total of spectra has an extra simple form,
in some sense.

In Section~2, we recall the expression of $\zeta_Q(2)$ 
given in \cite{IW},
and derive more explicit formula
of the generating function appearing in that expression.
We prove in Section~3 
our main results, the equations (\ref{eq:main1}) and (2).
The proof depends on several formulae
of hypergeometric series
 not only for ${}_2 F_1$
but also for ${}_3 F_2$
such as Clausen's identity.

\section{An expression of the generating function}

We start from the series-expression of the
special value
$\zeta_Q(2)$ of the non-commutative harmonic oscillator 
given in \cite[(4.5a)]{IW} 
\[
\zeta_Q(2) = Z_1(2) + \sum_{n=0}^\infty Z'_n(2).
\]
We introduce notations.
Recall that $\alpha>0$, $\beta>0$ with $\alpha\beta>1$.
Let us introduce the parameters
$\gamma = 1/\sqrt{\alpha\beta}$ and
$a=\gamma/\sqrt{1-\gamma^2}=1/\sqrt{\alpha\beta-1}$
as in \cite[(4.1)]{IW}.
Note that they satisfy $0<\gamma<1$ and $a>0$.

The term $Z_1(2)$ is given in \cite[(4.5b)]{IW} 
and $Z'_n(2)$ are given in \cite[(4.9)]{IW}
as
\begin{eqnarray}
Z_1(2) &=& \frac{(\alpha^{-1}+\beta^{-1})^2}{2(1-\gamma^2)}3 \zeta(2), \\
Z'_n(2) &=& (-1)^n \frac{(\alpha^{-1}-\beta^{-1})^2}{(1-\gamma^2)}
\binom{2n-1}{n} \left( \frac a2 \right)^{2n} J_n. 
\end{eqnarray}
The values $\{ J_n \}_{n=1,2,\cdots}$ are
specified by the generating function 
\[
w(z) :=\sum_{n=0}^\infty J_n z^n.
\]
The function $w(z)$
is a solution 
of the ordinary differential equation
\begin{equation}\label{eq:1}
z(1-z)^2 \frac{d^2 w}{dz^2} +(1-3z) (1-z) \frac{dw}{dz}
+\left(z-\frac34 \right) w = 0
\end{equation}
which is
given in \cite[Theorem~4.13]{IW}
and called a singly confluent Heun's differential equation.
The constant term is given by $w(0)=J_0 = 3 \zeta(2)=\pi^2/2$.
It is easy to see that there exists
a unique power-series solution of this homogeneous differential equation 
(\ref{eq:1})
with the initial condition $w(0)=\pi^2/2$.
The final target $\zeta_Q(2)$
involves these $J_n$'s with an infinite sum 
which does not seem
to have a closed form.

In this section,
we give a simple expression of the generating function $w(z)$.
We denote by $\pz=\partial/\partial z$.
\begin{lemma}
The differential equation (\ref{eq:1}) 
is equivalent to
\begin{equation}\label{eq:2}
4(1-z)\pz z \pz (1-z) w + w=0.
\end{equation}
\end{lemma}
Proof: This directly follows from Leibniz rule.
QED

\begin{lemma}
Let
$t=z/(z-1)$
be a new independent variable,
and  $\eta(t)=(1-z) w(z)$ a new unknown function.
Then the differential equation (\ref{eq:2})
is equivalent to
\begin{equation}\label{eq:3}
t(1-t) \pt^2 \eta + (1-2t) \pt \eta - \frac 14 \eta = 0.
\end{equation}
\end{lemma}
Proof:
The differential equation (\ref{eq:2}) is equivalent to
\[
4(z-1)^2\pz z \pz (z-1) w + (z-1)w=0.
\]
Note that $(z-1)(t-1)=1$ and $\pt:=\partial/\partial t=-(z-1)^2 \pz$.
Then
\[
4 \pt t(t-1) \pt \eta+ \eta=0.
\]
By Leibniz rule, this is equivalent to (\ref{eq:3}).
QED

\begin{proposition}\label{prop:1}
\[
w(z)  = \frac{J_0}{1-z} {} _2F_1\left(\frac12, \frac12; 1; \frac{z}{z-1} \right).
\]
\end{proposition}
Proof: Since any power-series solution of (\ref{eq:3}) in $t$
is a constant multiple of ${}_2F_1(\frac12,\frac12;1;t)$,
we have the conclusion. 
QED

\section{The special value}
We introduce the auxiliary series
\[
g(a) :=
\frac{2}{J_0} \sum_{n=0}^\infty (-1)^n
\binom{2n-1}{n} \left( \frac a2 \right)^{2n} J_n
\]
so that
\begin{eqnarray}
\zeta_Q(2) 
&=& \frac{(\alpha^{-1}+\beta^{-1})^2}{2(1-\gamma^2)} 3 \zeta(2)
+ \frac{(\alpha^{-1}-\beta^{-1})^2}{2(1-\gamma^2)} 3 \zeta(2) g(a) \\
&=& \frac{\pi^2}{4} \frac{(\alpha^{-1}+\beta^{-1})^2}{(1-\alpha^{-1}\beta^{-1})}
\left( 1+ \left(\frac{\alpha^{-1}-\beta^{-1}}{\alpha^{-1}+\beta^{-1}}\right)^2 g(a) \right)
\label{eq:zeta}
\end{eqnarray}

\begin{theorem}
\label{theorem}
\[
g(a) = 
{}_2F_1\left(\frac14,\frac34;1;- a^2\right)^2.
\]
\end{theorem}
Proof:
We note that
\[
\binom{2n-1}{n} \left(\frac{1}{2}\right)^{2n}
= \frac12 \times \frac{(2n-1)!!}{(2n)!!}
= \frac{1}{2\pi}\int_0^1 \frac{u^{n} du}{\sqrt{u(1-u)}}.
\]
%
%
Then,  integration by parts implies that
\begin{eqnarray}\label{eq:g}
g(a) &=& \frac{2}{2\pi J_0} \sum_{n=0}^\infty (-1)^n
\int_0^1 \frac{u^{n} du}{\sqrt{u(1-u)}}
 a^{2n} J_n 
= \frac{1}{\pi J_0} \int_0^1  
\frac{w(-a^2 u) du}{\sqrt{u(1-u)}}.
\end{eqnarray}
By Proposition~\ref{prop:1},
the function $w$ is written in terms of hypergeometric series $_2F_1$.
Substitute such an expression into the equation (\ref{eq:g}),
to obtain
\[
g(a) =  \frac{1}{\pi}
\int_0^1 
\frac{1}{1+a^2 u} {} _2F_1\left(\frac12, \frac12; 1; \frac{a^2u}{a^2 u+1} \right)
\frac{du}{\sqrt{u(1-u)}}.
\]
We introduce a new variable $v=(1+a^2)u/(1+a^2u)$.
Then
\[
g(a) = \frac{1}{\pi} \int_0^1 
{} _2F_1\left(\frac12, \frac12; 1; \frac{a^2v}{1+a^2} \right)
\frac{dv}{\sqrt{v(1-v)(1+a^2)}}.
\]
Now we use the formula (2.2.2) of \cite{AAR}
\[
{}_3F_2(a_1,a_2,a_3; b_1,b_2; x)
= \frac{\Gamma(b_2)}{\Gamma(a_3)\Gamma(b_2-a_3)}
\int_0^1 t^{a_3-1}(1-t)^{b_2-a_3-1} {}_2F_1(a_1,a_2; b_1; xt) dt.
\]
This shows that
\[
g(a) = \frac{1}{\sqrt{1+a^2}}
{}_3F_2 \left(\frac12,\frac12,\frac12; 1,1; \frac{a^2}{1+a^2} \right).
\]
By Clausen's identity (in e.g., Exercise~13 of Chapter 2 in \cite{AAR})
\[
{}_2F_1\left(a,b;a+b+\frac12; x\right)^2 
= {}_3F_2\left(2a,2b,a+b; 2a+2b,a+b+\frac12; x\right),
\]
we find that
\[
g(a) = \frac{1}{\sqrt{1+a^2}}
{}_2F_1\left(\frac14,\frac14;1;\frac{a^2}{1+a^2} \right)^2.
\]
Moreover
Pfaff  formula's Theorem 2.2.5 of \cite{AAR}
\[
{}_2F_1(a,b;c;x) = (1-x)^{-a} {}_2F_1(a,c-b;c;x/(x-1)),
\]
yields
\[
{}_2F_1 \left(\frac14,\frac34;1;-a^2\right) =
(1+a^2)^{-1/4} {}_2F_1\left(\frac14, \frac14;1; \frac{a^2}{a^2+1} \right).
\]
This shows that
\[
g(a) = 
{}_2F_1 \left( \frac14,\frac34;1;-a^2 \right)^2.
\]
QED


\begin{remark}
In the earlier version of the paper,
it was suggested to
make use of the hypergeometric series ${}_3 F_2$
with this special parameter $(1/2,1/2,1/2;1,1)$
by the multi-variable hypergeometric function of type $(3,6)$, 
especially by its restriction on the stratum called $X_{1b}$ in \cite{MSY}.
However,
we can avoid to use a multi-variable hypergeometric function
in the present version as is seen above.
\end{remark}

Theorem~\ref{theorem} 
with the help of the equation (\ref{eq:zeta})
shows the equation (\ref{eq:main1}).
The equation (\ref{eq:main2}) is shown as follows.
By Theorem 3.13 of \cite{AAR}
\[
{}_2F_1(a,b;2a;x) 
= \left(1-\frac{x}{2}\right)^{-b} 
{}_2F_1\left(\frac{b}{2},\frac{b+1}{2}; a+\frac{1}{2}; \left(\frac{x}{2-x}\right)^2\right),
\]
we have
\[
{}_2F_1\left(\frac12,\frac12;1; \frac{2ia}{ia+1}\right)
=(1+ia)^{1/2} {}_2F_1\left(\frac14,\frac34;1;-a^2\right).
\]
Let us recall the definition of the elliptic integral of the first kind;
\[
K(k) = \int_0^{\pi/2} \frac{d\theta}{\sqrt{1-k^2\sin^2\theta}}
= \frac{\pi}{2} {}_2F_1\left(\frac12,\frac12;1;k^2\right).
\]
Then we have
\[
{}_2F_1 \left( \frac14,\frac34;1;-a^2 \right)
=\frac{2}{\pi} (1+ia)^{-1/2} K\left( \frac{2ia}{ia+1} \right)
=\frac{2}{\pi} \int_0^{\pi/2} \frac{d\theta}{\sqrt{1+ia \cos2\theta}}
=\frac{1}{2\pi} \int_0^{2\pi} \frac{d\theta}{\sqrt{1+ia \cos\theta}},
\]
and the equation (\ref{eq:main2}).



\noindent
Deparment of Mathematics, Nagoya University\\
Chikusa, Nagoya 464-8602, Japan.\\
E-mail: ochiai@math.nagoya-u.ac.jp


\begin{thebibliography}{xxx}
\bibitem{AAR}
G. E. Andrews, R. Askey and R. Roy, 
 Special functions. Encyclopedia of Mathematics and its Applications, 71. Cambridge University Press, Cambridge, 1999.


\bibitem{IW}
T. Ichinose and M.Wakayama, 
Special values of the spectral zeta function
of the non-commutative harmonic oscillator and
confluent Heun equations,
Kyushu J. Math. {\bf 59} (2005), no. 1, 39--100.

\bibitem{MSY}
K. Matsumoto, T. Sasaki and M. Yoshida,  
The monodromy of the period map of a $4$-parameter family of $K3$ surfaces and the hypergeometric function of type $(3,6)$,
Internat. J. Math. {\bf 3} (1992), no. 1, 164 pp.

\bibitem{O}
H. Ochiai,
Non-commutative harmonic oscillators 
and Fuchsian ordinary differential operators,
Comm. Math. Phys. {\bf 217} (2001), 357--373.

\bibitem{O2}
H. Ochiai,
Non-commutative harmonic oscillators and
the connection problem for the Heun differential equation,
Letter in Math. Phys. {\bf 70} (2004), 133--139.

\bibitem{PW0}
A. Parmeggiani and M. Wakayama, 
Oscillator representations and systems of ordinary differential equations,
Proc. Nat. Acad. Sci. U.S.A. {\bf 98}(2001), 26--30.

\bibitem{PW1}
A. Parmeggiani and M. Wakayama,
Non-commutative harmonic oscillators. I, II.
Forum Math. {\bf 14} (2002), 539--604, 669--690.



\end{thebibliography}
\end{document}